\input amstex
\documentstyle{amsppt}

\input label.def
\input degt.def

\input epsf
\def\picture#1{\epsffile{#1-bb.eps}}
\def\cpic#1{$\vcenter{\hbox{\picture{#1}}}$}

\def\ie{\emph{i.e.}}
\def\eg{\emph{e.g.}}
\def\cf{\emph{cf}}
\def\via{\emph{via}}
\def\etc{\emph{etc}}

{\catcode`\@11
\gdef\proclaimfont@{\sl}}

\newtheorem{Remark}\Remark\thm\endAmSdef
\conjecture\thm\endproclaim

\loadbold
\def\bA{\bold A}
\def\bD{\bold D}
\def\bE{\bold E}
\def\bJ{\bold J}

\def\bY{\bold Y}
\def\bW{\bold W}

\let\Gl\lambda
\let\Go\omega
\let\Ge\epsilon
\let\Gk\kappa
\def\bGk{\bar\Gk}
\let\Gt\tau
\def\tGt{\tilde\Gt}
\let\tGt\Gt

\let\splus\oplus
\let\bigsplus\bigoplus
\let\ocup\cup
\let\bigocup\bigcup

\def\F{\Bbb F}

\def\CG#1{\Z_{#1}}
\def\DG#1{\Bbb D_{#1}}
\def\SGSet{\Bbb S}
\def\SG#1{\SGSet_{#1}}
\def\AG#1{\Bbb A_{#1}}

\def\O{\operatorname{\text{\sl O\/}}}
\def\GL{\operatorname{\text{\sl GL\/}}}
\def\PGL{\operatorname{\text{\sl PGL\/}}}

\def\CK{\Cal K}

\def\CS{\Cal S}
\def\CH{\Cal H}
\def\CC{\Cal C}
\def\CM{\Cal M}
\def\tc{\tilde c}
\def\tX{\tilde X}
\def\tS{\tilde S}
\def\tSigma{\tilde\Sigma}

\def\B{\bar B}
\def\L{\bar L}

\def\PP{\Bbb P}
\def\Cp#1{\PP^{#1}}
\def\term#1-{$\DG{#1}$-}

\def\ls|#1|{\mathopen|#1\mathclose|}
\let\<\langle
\let\>\rangle
\def\discr{\operatorname{discr}}
\def\Aut{\operatorname{Aut}}
\def\Pic{\operatorname{Pic}}
\def\Sym{\operatorname{Sym}}
\def\StSym{\Sym_{\roman{st}\!}}
\def\Sk#1{\operatorname{Sk}#1}

\topmatter

\author
Alex Degtyarev
\endauthor

\title
Stable symmetries of plane sextics
\endtitle

\address
Department of Mathematics,
Bilkent University,
06800 Ankara, Turkey
\endaddress

\email
degt\@fen.bilkent.edu.tr
\endemail

\abstract
We classify projective symmetries of irreducible plane sextics
with simple singularities which are stable under equivariant
deformations. We also outline a connection between order~$2$
stable symmetries and maximal trigonal curves.
\endabstract



\endtopmatter

\document

\section{Introduction}

\subsection{Motivation}
In a recent series of papers~\cite{degt.Oka3}, \cite{degt.8a2},
\cite{degt-Oka}, we described the moduli spaces and
computed the fundamental groups of a number of plane sextics. A
common feature of all these papers is the fact that we start with
proving that \emph{each} sextic in the equisingular stratum under
consideration possesses a certain symmetry (projective
automorphism); then, this symmetry is used to analyze the moduli
space and to write down explicit equations defining each curve.
(Note that, prior to the papers cited above, for most curves in
question only
the existence was known, which was proved by rather indirect
means.) In most cases, the symmetry can also be used to facilitate
the computation of the fundamental group: in a certain sense, one
can reduce pencils of degree six to those of degree four.

Thus, a natural question arises: {\it are there other plane
sextic curves admitting a symmetry stable under equisingular
deformations?} In the present paper, we give a partial answer to
this question: our principal result is Theorem~\ref{th.list},
classifying all \emph{irreducible} plane sextics with simple
singularities whose group of stable symmetries is nontrivial.
Our approach should also
apply to reducible sextics with simple singularities,
but the number of cases to be considered and the number of
classes obtained should be much larger. On the other hand, it
appears that sextics with a non-simple singular point do not admit
stable symmetries, \cf. Remark~\ref{non-simple}.

In a subsequent paper, we are planning to use the results obtained
here to compute the fundamental groups of all sextics listed in
Theorem~\ref{th.list}.

\subsection{Stable symmetries}
The (combinatorial) set of singularities of a plane
curve $B$ is denoted by $\Sigma=\Sigma_B$. If all singularities
are simple, we identify~$\Sigma$ with the lattice spanned
by the classes of exceptional divisors in the minimal resolution
of singularities of the double plane ramified at~$B$,
\cf~\ref{s.K3}; for this reason, we use $\splus$ in
the notation.

For an integer $m\ge1$, we denote by $\CC_m$ the space of plane curves
of degree~$m$; as is well known, $\CC_m$ is a projective space of
dimension $n(n+3)/2$. Given a set of
singularities~$\Sigma$, the \emph{equisingular stratum}
$\CC_m(\Sigma)\subset\CC_m$ is the set of curves whose
set of singularities is~$\Sigma$.
(Throughout this paper, all singularities are isolated and, in
most cases, even simple.)
There is an
obvious action of $\PGL(3,\C)$ on~$\CC_m$ preserving each stratum
$\CC_m(\Sigma)$; the quotient
$\CM_m(\Sigma)=\CC_m(\Sigma)/\!\PGL(3,\C)$ is called the
\emph{moduli space} of curves of degree~$m$ with the set of
singularities~$\Sigma$.

\definition\label{def.symmetry}
A \emph{symmetry} of a plane curve $B\subset\Cp2$ is an
automorphism of the pair $(\Cp2,B)$. (We use the term
`symmetry' instead of `automorphism' to
avoid confusion with automorphisms of~$B$ as abstract curve.)
A symmetry
$s$ of~$B$ is called \emph{stable} if there is a neighborhood $U$
of~$B$ in the equisingular stratum $\CC_m(\Sigma_B)$
and a continuous family
$\sigma\:U\to\PGL(3,\C)$ such that $\sigma(B')$ is a symmetry
of~$B'$ for each $B'\in U$ and $\sigma(B)=s$.
The set of (stable) symmetries of~$B$ is denoted by $\Sym B$
(respectively, $\StSym B$).
\enddefinition

Alternatively, one can consider the bundle
$$
\Sym(\Sigma)=\bigl\{(B,s)\in\CC_m(\Sigma)\times\PGL(3,\C)\bigm|
 s\in\Sym B\bigr\}\to\CC_m(\Sigma)
$$
and define the (nonabelian) sheaf $\StSym(\Sigma)$ of germs of
sections of $\Sym(\Sigma)$. Then $\StSym B$ is the stalk of
$\StSym(\Sigma)$ over~$B$. Note that \emph{a priori} it is not
obvious that $\StSym(\Sigma)$ is locally constant or even that the
groups $\StSym B$ are semicontinuous in any reasonable sense.
Below we show, see Corollary~\ref{locally.constant},
that the sheaf $\StSym(\Sigma)$ is indeed locally
constant in the case of curves of degree six with simple
singularities only. Furthermore, we show that, in this case, the
group of stable symmetries is the group of symmetries of a generic
curve in a given equisingular stratum (see
Corollary~\ref{Sym=StSym}); thus, the study of $\StSym B$ is
equivalent to the study of a generic fiber $\PGL(3,\C)/\StSym B$
of the projection $\CC_6(\Sigma_B)\to\CM_6(\Sigma_B)$.

It is worth mentioning that stable symmetries can also be used to
describe the moduli space $\CM_6(\Sigma_B)$ itself. Thus, from the
results of~\cite{degt.Oka3}, \cite{degt.8a2},
and~\cite{degt-Oka} it follows that, for all sets of singularities
mentioned in Theorem~\iref{th.list}{w=9}, \ditto{14}, \ditto{w=8},
and~\ditto{10} below, the moduli spaces are unirational. It is
anticipated that a similar statement holds for (most of) the other
sets of singularities listed in Theorem~\ref{th.list}; we will
discuss this in details in a subsequent paper.

\subsection{Principal results}
In order to state our principal result, we introduce a few terms.
A \emph{\term2p-sextic} is an irreducible plane sextic~$B$ with simple
singularities and such that the fundamental group
$\pi_1(\Cp2\sminus B)$ factors
to the dihedral group~$\DG{2p}$. As shown in~\cite{degt.Oka},
there are \term6-, \term10-, and \term14-sextics. Furthermore, the
class of \term6-sextics coincides with the class of irreducible
sextics with simple singularities that are of \emph{torus type},
\ie, whose equation can be represented in the form $p^3+q^2=0$,
where $p$ and~$q$ are some homogeneous polynomials of degree~$2$
and~$3$, respectively.

There are relatively few \term10- and \term14-sextics: they form,
respectively, eight and two deformation families, see~\cite{degt.Oka}.
In order to
describe the hierarchy of \term6-sextics, we introduce the
\emph{weight} $w(B)$ as the total weight of the singularities
of~$B$, where the \emph{weight} of a simple singular point~$P$ is
defined as follows: $w(\bA_{3i-1})=i$, $w(\bE_6)=2$, and $w(P)=0$
otherwise. One has $w(B)\le9$, if $B$ is of torus type, then
$w(B)\ge6$, and a sextic of weight $\ge6$ is of torus type unless
$w(B)=6$ and all singular points of~$B$ of weight zero are nodes.
In the latter exceptional case, most sets of singularities are
realized by two deformation families, one of torus type and one
not; for the complete classification, see
A.~\"Ozg\"uner~\cite{Aysegul}.

If $B$ is a plane sextic of torus type and $p^3+q^2$ is a
\emph{torus structure} of~$B$, denote by $\Sigma_B^{\inj}$ the
set of inner singularities of~$B$. (Recall that a singular
point~$P$
of~$B$ is called \emph{inner} with respect to a given torus
structure $p^3+q^2$ if $P$ belongs to the intersection of the
conic $\{p=0\}$ and the cubic $\{q=0\}$.)

The principal result of the present paper
is the following theorem, which
gives a complete classification
of stable symmetries of irreducible plane sextics with simple
singularities.

\theorem\label{th.list}
The following is the complete list of irreducible
plane sextics with simple
singularities and nontrivial
group $\StSym B$ of stable symmetries\rom:
\roster
\item\local{w=9}
all sextics of weight nine\rom:
$\StSym B=(\CG3\times\CG3)\rtimes\CG2$,
the $\CG2$ factor
acting on the kernel $\CG3\times\CG3$
\via\ the multiplication by~$(-1)$\rom;
\item\local{3E}
\term6-sextics with
$\Sigma_B^{\inj}=3\bE_6$\rom:
$\StSym B$ is the symmetric group~$\SG3$\rom;
\item\local{14}
all \term14-sextics\rom: $\StSym B=\CG3$\rom.
\endroster
For the rest of the list, one has $\StSym B=\CG2$\rom:
\roster
\item[4]\local{E}
\term6-sextics with
$\Sigma_B^{\inj}=2\bE_6\splus\bA_5$ or $2\bE_6\splus2\bA_2$\rom;
\item\local{A}
\term6-sextics with
$\Sigma_B^{\inj}=\bA_{17}$ or $2\bA_8$\rom;
\item\local{w=8}
all sextics of weight eight\rom;
\item\local{10}
all \term10-sextics\rom;
\item\local{E8}
sextics with
$\Sigma_B=2\bE_8\splus\Sigma'$, where
$\Sigma'=\bA_3$, $\bA_2$, or $k\bA_1$, $k=0,1,2$.
\endroster
\endtheorem

Theorem~\ref{th.list} is proved in Section~\ref{proof.list}.

\Remark
In items~\loccit{3E}, \loccit{E}, and~\loccit{A},
the curves have weight~$\le7$;
hence, each curve~$B$ has a unique torus structure,
see~\cite{degt.Oka},
and $\Sigma_B^{\inj}$ is well defined. Item~\loccit{E8} lists
all sextics with two type~$\bE_8$ singular points, see
Proposition~\ref{2E8.classes}.
\endRemark

There is a mysterious connection between stable involutions of
irreducible plane sextics and maximal
(in the sense of~\cite{degt.Oka}, see
Definition~\ref{def.maximal} for details) trigonal curves in the
cone~$\Sigma_2$. Roughly, an involution $c\in\Sym B$ is stable if
and only if the quotient $B/c$ is maximal. We postpone the precise
statement till \S\ref{S.involutions}, see
Theorems~\ref{th.maximal} and~\ref{th.converse}, as they require a
number of preliminary definitions.

\subsection{Contents of the paper}
In \S\ref{S.reduction}, we reduce the problem of classification of
stable symmetries to a combinatorial question. First, we apply the
theory of $K3$-surfaces and describe the symmetries of plane
sextics with simple singularities in arithmetical terms, see
Theorem~\ref{th.Sym}. Next, in Theorem~\ref{stable=>symplectic},
we give an arithmetical characterization of stable symmetries.
With a few exceptions with the maximal total Milnor
number $\mu=19$, this theorem applies to reducible curves as well.
Finally, in Theorem~\ref{th.Dynkin}, we describe stable symmetries
of irreducible sextics in terms of symmetries of their Dynkin
graphs. With the few exceptions above, this theorem also applies
to reducible curves, provided that the definition of configuration
and its stable symmetry is modified to take into account the
hyperplane section class.

In \S\ref{S.proof}, we classify stable symmetries of Dynkin graphs
of irreducible sextics and prove Theorem~\ref{th.list}.

\S\ref{S.involutions} deals with stable involutions and trigonal
curves. First, we classify all
stable maximal trigonal curves in the
Hirzebruch surface~$\Sigma_2$, see
Theorem~\ref{th.stable.curves}. Then, comparing this result and
Theorem~\ref{th.list}, we give a characterization of
stable involutions of
irreducible plane sextics in terms of the maximality of the
quotient curve. There is strong evidence that a similar relation
holds as well for reducible sextics with simple singularities, see
Conjecture~\ref{conjecture} and Remark~\ref{conjecture.proof}.

\section{The combinatorial reduction}\label{S.reduction}

The principal result of this section is Theorem~\ref{th.Dynkin},
reducing the study of stable symmetries of plane sextics to the
study of symmetries of their Dynkin diagrams.

\subsection{Discriminant forms}\label{s.discr}
An \emph{\rom(integral\rom) lattice} is a finitely generated free
abelian group~$S$ supplied with a symmetric bilinear form
$b\:S\otimes S\to\Z$. We abbreviate $b(x,y)=x\cdot y$ and
$b(x,x)=x^2$. A lattice~$S$ is called \emph{even} if $x^2=0\bmod2$ for
all $x\in S$. As the transition matrix between two integral bases
has determinant $\pm1$, the determinant $\det S\in\Z$ (\ie, the
determinant of the Gram matrix of~$b$ in any basis of~$S$) is well
defined. A lattice~$S$ is called \emph{nondegenerate} if
$\det S\ne0$; it is called \emph{unimodular} if $\det S=\pm1$.

Given a lattice~$S$,
the form~$b$ extends to a form
$(S\otimes\Q)\otimes(S\otimes\Q)\to\Q$.
If
$S$ is nondegenerate, the dual group $S^*=\Hom(S,\Z)$ can
be identified with the subgroup
$$
\bigl\{x\in S\otimes\Q\bigm|
 \text{$x\cdot y\in\Z$ for all $x\in S$}\bigr\}.
$$
In particular, $S\subset S^*$
is a finite index subgroup. The quotient $S^*\!/S$
is called the \emph{discriminant group} of~$S$ and is
denoted by $\discr S$ or~$\CS$. The discriminant group
inherits from $S\otimes\Q$ a symmetric bilinear form
$b_{\CS}\:\CS\otimes\CS\to\Q/\Z$,
called the \emph{discriminant form}, and,
if $S$ is even, its quadratic extension~$q_{\CS}$, \ie, a function
$q_{\CS}\:\CS\to\Q/2\Z$ such that
$q_{\CS}(x+y)=q_{\CS}(x)+q_{\CS}(y)+2b_{\CS}(x,y)$
for all $x,y\in\CS$, where $2$
is regarded as a homomorphism $\Q/\Z\to\Q/2\Z$.
One has $\ls|\CS|=\ls|\det S|$; in particular,
$\CS=0$ if and only if $S$ is unimodular.

Given a prime~$p$, we use the notation $\discr_pS=\CS_p$ for the
$p$-primary part of~$\CS$. One has $\CS_p=\CS\otimes\Z_{p^r}$,
$r\gg1$, and $\CS=\bigoplus_p\CS_p$, the sum running over all
primes.

From now on, \emph{all lattices considered are even and
nondegenerate}.

An \emph{extension} of a lattice~$S$ is another lattice~$M$
containing~$S$, so that the form on~$S$ is the restriction of that
on~$M$. An \emph{isomorphism} between two extensions
$M_1\supset S$ and $M_2\supset S$ is an isometry $M_1\to M_2$
whose restriction to~$S$ is the identity.
Next three theorems are found in Nikulin~\cite{Nikulin}.

\theorem\label{th.Nik1}
Given a 
lattice~$S$, there is a canonical one-to-one correspondence
between the set of isomorphism classes of finite index extensions
$M\supset S$
and the set of isotropic subgroups $\CK\subset\CS$.
Under this correspondence, one has
$M=\{x\in S^*\,|\,x\bmod S\in\CK\}$ and
$\discr M=\CK^\perp\!/\CK$.
\qed
\endtheorem

The isotropic subgroup $\CK\subset\CS$ as in Theorem~\ref{th.Nik1}
is called the \emph{kernel} of the extension $M\supset S$. It can
be defined as the image of $M\!/S$ under the homomorphism induced by
the natural inclusion $M\hookrightarrow S^*$.

\theorem\label{th.Nik2}
Let $M\supset S$ be a finite index
extension of a
lattice~$S$,
and let $\CK\subset\CS$ be its
kernel.
Then, an auto-isometry $S\to S$ extends to~$M$ if and only if the
induced automorphism of~$\CS$ preserves~$\CK$.
\qed
\endtheorem

\theorem\label{th.Nik3}
Let $S\subset M$ be a primitive
sublattice of a unimodular
lattice~$M$. Then the kernel of the finite index extension
$M\supset S\oplus S^\perp$ is the graph of
an anti-isometry
$\discr S\to\discr S^\perp$.
\qed
\endtheorem

We will use Theorems~\ref{th.Nik2} and~\ref{th.Nik3}
in the following form.

\corollary\label{extension}
Let $S\subset M$ be a
sublattice of a unimodular
lattice~$M$, and let $\CK\subset\CS$ be the kernel of the
extension $\tilde S\supset S$, where $\tilde S$ is the primitive
hull of~$S$ in~$M$. Consider an auto-isometry $c\:S\to S$. Then,
$c\oplus\id_{S^\perp}$ extends to~$M$ if and only if $c$
preserves~$\CK$ and the auto-isometry of
$\CK^\perp\!/\CK$ induced by~$c$ is trivial.
\endcorollary

\proof
Apply Theorem~\ref{th.Nik2} twice: first, to
$\tilde S\supset S$, then to
$M\supset\tilde S\oplus S^\perp$.
\endproof

\subsection{Root systems}\label{s.root}
A \emph{root} in an even lattice~$S$ is an element $v\in S$ of
square~$-2$.
A \emph{root system} is a negative definite lattice
generated by its roots. Every root system admits a unique decomposition
into an orthogonal sum of irreducible root systems, the latter
being either $\bA_p$, $p\ge1$, or~$\bD_q$, $q\ge4$,
or~$\bE_6$, $\bE_7$, $\bE_8$. The discriminant forms are as
follows:
$$
\gathered
\discr\bA_p=[-\tfrac{p}{p+1}],\quad
\discr\bD_{2k+1}=[-\tfrac{2k+1}4],\\
\discr\bD_{8k\pm2}=2[\mp\tfrac12],\quad
\discr\bD_{8k}=\Cal U_2,\quad
\discr\bD_{8k+4}=\Cal V_2,\quad\\
\discr\bE_6=[\tfrac23],\quad
\discr\bE_7=[\tfrac12],\quad
\discr\bE_8=0.
\endgathered\eqtag\label{eq.discr}
$$
Here, $[\frac pq]$ is the cyclic group $\Z_q$ generated by
an element of square $\frac pq\in\Q/2\Z$, and $\Cal U_2$
(respectively, $\Cal V_2$) is the quadratic form on
$\Z_2\oplus\Z_2$ generated by elements~$x$, $y$ with
$x\cdot y=\frac12\in\Q/\Z$ and $x^2=y^2=0\in\Q/2\Z$
(respectively, $x^2=y^2=1$).

Given a root system~$S$, the group generated by reflections
(defined by the roots of~$S$) acts simply transitively on the set
of Weyl chambers of~$S$. The roots defining the walls of any
Weyl chamber form a \emph{standard basis} for~$S$. The incidence
graph~$\Gamma$ of a standard basis is called the \emph{Dynkin
diagram} of~$S$. Irreducible root systems correspond to connected
Dynkin diagrams. With a certain abuse of the language, we will
speak about the \emph{discriminant group} $\discr\Gamma$ of a Dynkin
diagram~$\Gamma$, referring to the discriminant group of the root
system~$S_\Gamma$ spanned by~$\Gamma$.

Denote by $\Sym\Gamma$ the group of symmetries of a Dynkin
diagram~$\Gamma$. There is an obvious homomorphism
$$
\discr\:\Sym\Gamma@>>>\O(S_\Gamma)@>>>\Aut\discr\Gamma.
$$
The following three statements are well known; they follow
immediately from the classification of connected Dynkin diagrams,
see, \eg, N\.~Bourbaki~\cite{Bourbaki}.

\lemma\label{discr.Dynkin}
Let $\Gamma$ be a connected Dynkin diagram. Then\rom:
\roster
\item
$\discr\Gamma\ne0$ unless $\Gamma$ is of type~$\bE_8$\rom;
\item
the homomorphism $\discr\:\Sym\Gamma\to\Aut\discr\Gamma$ is monic.
\qed
\endroster
\endlemma

\lemma\label{Sym.Dynkin}
Let $\Gamma$ be a connected Dynkin diagram. Then\rom:
\roster
\item
if $\Gamma$ is of type~$\bD_4$, then
$\Sym\Gamma=\Aut\discr\Gamma=\SG3$\rom;
\item
if $\Gamma$ is of type~$\bA_1$, $\bE_7$, or~$\bE_8$, then
$\Sym\Gamma=\Aut\discr\Gamma=1$\rom;
\item
for all other types, $\Sym\Gamma=\CG2$.
\qed
\endroster
\endlemma

\lemma\label{-id.Dynkin}
If $\Gamma$ is a connected Dynkin diagram
of type~$\bA_p$, $p\ge2$, $\bD_{2k+1}$, or~$\bE_6$,
then the only nontrivial symmetry of~$\Gamma$ induces $-\id$ on
$\discr\Gamma$.
\qed
\endlemma

Further details on irreducible root systems are found in
Bourbaki~\cite{Bourbaki}.

\subsection{The covering $K3$-surface}\label{s.K3}
Let $B$ be a plane sextic with simple singularities.
Denote by $X_B\to\Cp2$ the double covering ramified at~$B$, and
let~$\tX_B$ be the minimal resolution of~$X_B$. It is a
$K3$-surface. Let $\tGt\:\tX\to\tX$ be the deck translation
of the covering $\tX_B\to\Cp2$, and let $\Aut_{\tGt}\tX_B$ be the
centralizer of~$\tGt$ in the group of automorphisms of~$\tX_B$.
Since any symmetry of~$B$ lifts to two automorphisms of~$\tX$,
there is an exact sequence
$$
1@>>>\{\id,\tGt\}@>>>\Aut_{\tGt}\tX_B@>>>\Sym B@>>>1.
\eqtag\label{eq.AutX}
$$

\definition
The \emph{homological type} of a sextic~$B$ with simple
singularities is the triple $\CH_B=(L_B,h_B,\Gamma_B)$, where
$L_B$ is the lattice $H_2(\tX_B)$, $h_B\in L_B$ is the class of
the pull-back of a generic line (so that $h^2=2$), and $\Gamma_B$
is the set of classes of the exceptional divisors over the
singular points of~$B$.

An \emph{automorphism} of the homological type
$\CH_B=(L_B,h_B,\Gamma_B)$ is an isometry of~$L_B$
preserving~$h_B$ and~$\Gamma_B$ (as a set). The group of
automorphisms of~$\CH_B$ is denoted by $\Aut\CH_B$. We will also
consider the subgroups
$\Aut_+\CH_B$ and $\Aut_\pm\CH_B$ consisting of
the automorphisms inducing, respectively,~$\id$ and $\pm\id$ on the
orthogonal complement $(\Gamma_B\cup h_B)^\perp$.
\enddefinition

Denote by $\Sigma_B\subset L_B$ the sublattice
spanned by~$\Gamma_B$, and let $S_B=\Sigma_B\oplus\<h_B\>$. The
primitive hulls of~$\Sigma_B$ and~$S_B$ in~$L_B$ are denoted
by~$\tSigma_B$ and~$\tS_B$, respectively, and the kernel of the
finite index extension $\tS_B\supset S_B$ is denoted by~$\CK_B$.

\lemma\label{Aut+-}
$\Aut_\pm\CH_B$ is the subgroup of elements $s\in\Aut\CH_B$
inducing a scalar on $S_B^\perp$.
\endlemma

\proof
For any $s\in\Aut\CH_B$, both the restriction $s|_{S^\perp}$ and
its inverse are defined over~$\Z$; hence, if $s|_{S^\perp}$ is a
scalar, one must have $s|_{S^\perp}=\pm\id$.
\endproof

The following statement is contained in~\cite{JAG}.

\proposition\label{irreducible}
If $B$ is irreducible, then $\CK_B$ is free of $2$-torsion. In
particular, one has $\CK_B\subset\discr\Sigma_B$,
$\tS_B=\tSigma_B\oplus\<h_B\>$, and
$\discr\tS_B=\discr\tSigma_B\oplus\<\frac12h_B\>$.
\qed
\endproposition

The lattice~$\Sigma_B$ is a root system, and the elements
of~$\Gamma_B$ form a standard basis for~$\Sigma_B$,
see~\ref{s.root}.
In what follows, we identify the set~$\Gamma_B$ with its incidence
graph. The connected components of~$\Gamma_B$ (irreducible
components of~$\Sigma_B$) are in a one-to-one correspondence with
the singular points of~$B$, each component being a connected
Dynkin diagram (respectively, irreducible root
system) of the same name as the type of the singular
point.

\definition
The \emph{period} of a sextic~$B$ is the $1$-subspace
$\Go_B=H^{2,0}(\tX_B)\subset L_B\otimes\C$; it is formed by the
classes of holomorphic $2$-forms on~$\tX_B$. The
\emph{extended homological type} of~$B$ is the pair
$(\CH_B,\Go_B)$. An \emph{automorphism} of the extended
homological type is an automorphism of~$\CH_B$
preserving~$\Go_B$. The group of automorphisms of
$(\CH_B,\Go_B)$ is denoted by $\Aut(\CH_B,\Go_B)$.
\enddefinition

Recall that the period~$\Go_B$ is a point in the projectivization
of the cone
$$
\Omega=\{x\in S_B^\perp\otimes\C\,|\,x^2=0,\,x\cdot\bar x>0\}.
\eqtag\label{eq.cone}
$$
Conversely, any generic (complementary to a countable union of
hyperplanes) point in the projectivization of~$\Omega$ is the
period of a certain sextic, which is in the same equisingular
stratum $\CC_6(\Sigma_B)$ as~$B$.

There are obvious inclusions
$$
\Aut_+\CH_B\subset\Aut_\pm\CH_B\subset\Aut(\CH_B,\Go_B)\subset\Aut\CH_B.
$$
By definition, $\Go_B$ is an eigenspace of any element
$a\in\Aut(\CH_B,\Go_B)$. Sending~$a$ to the corresponding
eigenvalue defines a homomorphism
$\Gl\:\Aut(\CH_B,\Go_B)\to\C^*$.

Any automorphism $\tc\in\Aut_{\tGt}\tX_B$
induces an automorphism
$\tc_*\in\Aut(\CH_B,\Go_B)$. In particular, $\tGt$ itself
induces an automorphism
$\tGt_*\in\Aut(\CH_B,\Go_B)$, which can be described as
follows. On each connected component of~$\Gamma_B$ of
type~$\bA_p$, $p\ge2$, $\bD_{2k+1}$, or~$\bE_6$,
$\tGt_*$ is the only nontrivial symmetry, on any other component
and on~$\<h_B\>$, it is the identity, and
on~$S_B^\perp$, minus identity. The map just described
preserves~$\CK_B$ (as well as any subgroup of $\discr S_B$), and
the induced automorphisms of all discriminants are $-\id$;
due to Theorem~\ref{th.Nik2},
the map extends to~$L_B$.

\theorem\label{th.Sym}
For any plane sextic $B\in\Cp2$ with simple singularities,
the map $\tc\mapsto\tc_*$ establishes an
isomorphism $\Aut_{\tGt}\tX_B=\Aut(\CH_B,\Go_B)$. Hence,
there is an exact sequence
$$
1@>>>\{\id,\tGt_*\}@>>>\Aut(\CH_B,\Go_B)@>>>\Sym B@>>>1,
$$
obtained from~\eqref{eq.AutX} \via\ the above isomorphism.
\endtheorem

\proof
Let $\Pic\tX_B=\Go_B^\perp\cap L_B$ be the Picard group
of~$\tX_B$. Recall that the K\"ahler cone $V_B^+$ of~$\tX_B$ can
be
defined as the set
$$
\bigl\{x\in \Go_B^\perp\cap(L_B\otimes\R)\bigm|
\text{
$x^2>0$ and $x\cdot[E]>0$ for any $(-2)$-curve~$E\subset\tX_B$}
\bigr\}.
$$
The projectivization $\PP(V^+_B)$ is an (open) fundamental
polyhedron of the group of motions of the hyperbolic space
$\PP(\{x\in\Go_B^\perp\cap(L_B\otimes\R)\,|\,x^2>0\})$
generated by the reflections
defined by the roots of $\Pic\tX_B$. In the
case under consideration, $V_B^+$
is characterized
(among the other fundamental polyhedra)
by the following properties:
\roster
\item
$V_B^+\cdot v>0$ for any $v\in\Gamma_B$;
\item
the closure of~$V_B^+$ contains~$h_B$.
\endroster
Consider an element $\tc_*\in\Aut(\CH_B,\Go_B)$ and
regard it as an isometry of $H_2(\tX_B)$. By definition,
$\tc_*$ preserves~$h_B$, $\Gamma_B$, and~$\Go_B$; hence,
$\tc_*$ also preserves~$V_B^+$. Now, a standard argument using the
description of the
fine period space of marked K\"ahler $K3$-surfaces, see
A.~Beauville~\cite{Beauville}, shows that any isometry~$\tc_*$
of $H_2(\tX_B)$ preserving $\Go_B$ and $V_B^+$ is induced by a
unique automorphism~$\tc$ of~$\tX_B$. Since $\tc_*(h_B)=h_B$
and $h_B$ (regarded as an element of $\Pic\tX_B$) is the linear
system defining the projection $\tX\to\Cp2$, one has
$\tc\in\Aut_{\tGt}\tX_B$. The existence (uniqueness) of~$\tc$
above assert that the map $\tc\mapsto\tc_*$ is onto
(respectively, one-to-one).
\endproof

\proposition
For any sextic~$B$, the group $\Sym B$ is finite.
\endproposition

\proof
Since $\tX_B$ is obviously algebraic, the kernel of the canonical
representation $\Aut\tX_B\to\O(\Pic\tX_B)$ is a finite cyclic
group, see Nikulin~\cite{Nikulin.auto}. On the other hand, the
image of $\Aut_{\tGt}\tX_B$ is a subgroup of
$\O(h_B^\perp\cap\Pic\tX_B)$. Since
$h_B^\perp\cap\Pic\tX_B$ is a negative definite lattice,
its group of isometries is finite.
\endproof

\subsection{The symplectic lift}
Recall that an automorphism of a $K3$-surface~$X$ is called
symplectic (anti-symplectic) if it preserves (respectively,
reverses) holomorphic $2$-forms on~$X$. Note that any automorphism
multiplies all $2$-forms by a certain constant $\Gl\in\C^*$; the
automorphism is (anti-)symplectic if and only if $\Gl=\pm1$.

\theorem\label{stable=>symplectic}
Let $B\subset\Cp2$ be a sextic with simple singularities, and assume
that either $\mu(B)<19$ or $B$ is irreducible.
Then, for any stable symmetry~$c$ of~$B$,
one of the two lifts of~$c$ to the
covering $K3$-surface $\tX_B$ is symplectic, and the other one is
anti-symplectic.
The induced automorphisms of~$\CH_B$ belong to $\Aut_\pm\CH_B$.
\endtheorem

\proof
Since the two lifts differ by~$\tGt$,
which is anti-symplectic and induces $-\id$ on
$S_B^\perp$,
it suffices to show that any lift induces $\pm\id$ on $S_B^\perp$.
(Then it acts \via\ $(\pm1)$ on~$\Go_B$ and, hence, is symplectic or
anti-symplectic.)

Let $\tc_*$ be the automorphism of $S_B^\perp$ induced by the
chosen lift. Then the period $\Go_B\subset S_B^\perp\otimes\C$ is an
eigenspace of~$\tc_*$.
Assume that $\rank S_B^\perp\ge3$.
Since $c$ is stable,
there is a neighborhood~$U$ of~$\Go_B$ in the projectivization of
the cone~$\Omega$, see~\eqref{eq.cone},
such that any generic $\Go\in U$ is also an eigenspace of~$\tc_*$,
obviously corresponding to the same eigenvalue
as~$\Go_B$.
On the other
hand, any open subset of~$\Omega$ spans $S_B^\perp\otimes\C$.
Thus, $\tc_*$ is a scalar and Lemma~\ref{Aut+-} applies.

Now, assume that $\rank S_B^\perp=2$. Any positive definite
lattice of rank~$2$ admitting an orientation preserving
automorphism other than $\pm\id$ is isomorphic to
either
$\bA_2(-m)$ or
$2\bA_1(-m)$, where $m$ is a positive integer and $S(m)$ means
that the bilinear form on the lattice~$S$ is multiplied by~$m$.
(For example, one can argue that
these are the lattices in $\C^1=\R^2$ admitting
a non-trivial complex multiplication.)
On the other hand, since $B$ is irreducible,
the discriminant $\discr\smash{\tS_B}\cong-\discr S_B^\perp$
(see Theorem~\ref{th.Nik3}) has a direct summand
$\<\frac12h_B\>\cong[\frac12]$, see Proposition~\ref{irreducible}.
It is
immediate that neither $\discr\bA_2(-m)$ nor $\discr2\bA_1(-m)$,
$m>0$, has a direct summand $[-\frac12]$.
\endproof

\example
There do exist reducible plane sextics
for which the conclusion of
Theorem~\ref{stable=>symplectic} fails.
For example, take
for~$B$ the curve given, in some affine coordinates $(x,y)$,
by the equation
$$
(y^3-y)(y^3-y+x^3)=0.
$$
The set of singularities of~$B$ is $\bD_4\splus3\bA_5$; hence,
$\mu(B)=19$ and $\StSym B=\Sym B$.
(Note also that $B$ is a sextic of torus type with four
irreducible components; a simple calculation shows that
$S_B^\perp=\bA_2(-1)$.)
An affine part of~$X_B$ is given by
$$
z^2=(y^3-y)(y^3-y+x^3).
$$
The lift $(x,y,z)\mapsto(\epsilon x,y,z)$ of the symmetry
$(x,y)\mapsto(\epsilon x,y)$, $\epsilon^3=1$, is an order~$3$
automorphism of~$X_B$ with a dimension one component $\{x=0\}$
in the fixed
point set; hence, it is neither symplectic nor anti-symplectic.
\endexample

\proposition\label{symplectic=>stable}
For any automorphism $\tc_*\in\Aut_\pm\CH_B$, its image
in $\Sym B$,
see Theorem~\ref{th.Sym}, is stable.
\endproposition

\proof
For any sextic~$B'$ close to~$B$ in its equisingular
stratum $\CC_6(\Sigma_B)$, one can identify
$\CH_{B'}$ and $\CH_B$ and, under this identification, $\Go_{B'}$
is a $1$-space close to~$\Go_B$ in $\PP(S_B^\perp\otimes\C)$.
Hence, $\tc_*$ is also an automorphism of
$(\CH_{B'},\Go_{B'})$, and $c$ is stable.
\endproof

\corollary\label{Sym=StSym}
If $B$ is generic in its equisingular stratum
\rom(\ie, $B$
belongs to the complement of a certain countable union of
codimension~$1$ subsets of $\CC_m(\Sigma_B)$\,\rom),
then
$\StSym B=\Sym B$.
\endcorollary

\proof
If $\mu(B)=19$, the statement is obvious, as the moduli space
$\CM_6(\Sigma_B)$ is discrete. Otherwise,
in view of Theorem~\ref{th.Sym} and
Proposition~\ref{symplectic=>stable}, it suffices to show that
$\Aut(\CH_B,\Go)=\Aut_\pm\CH_B$ for a generic element
$\Go\in\PP(\Omega)$, see~\eqref{eq.cone}. The latter statement
follows from the the fact that the group $\Aut\CH_B$ is countable
and from Lemma~\ref{Aut+-}, which implies that, for any
$s\in\Aut\CH_B\sminus\Aut_\pm\CH_B$, the eigenspaces of the
restriction $s|_{S^\perp}$ are proper subspaces of~$S_B^\perp$.
\endproof

\subsection{Reduction to Dynkin diagrams}
From now on, we consider irreducible sextics only
and reserve the notation~$\tc$ for the \emph{symplectic}
lift of a stable symmetry~$c$ to the covering $K3$-surface~$\tX_B$;
it is well defined due to Theorem~\ref{stable=>symplectic}. We
denote by~$\tc_*$ the
induced isometry of $L_B=H_2(\tX_B)$, and by
$\tc_\Gamma\:\Gamma_B\to\Gamma_B$, the
induced symmetry of the Dynkin diagram.

\definition\label{def.configuration}
The \emph{configuration} of an irreducible plane
sextic~$B$ with simple
singularities is the pair $(\Gamma_B,\CK_B)$, where
$\CK_B\subset\discr\Gamma_B$ is the kernel of the extension
$\tSigma\supset\Sigma$, see~\ref{s.K3}. A \emph{symmetry} of
$(\Gamma_B,\CK_B)$ is a symmetry $s\in\Sym\Gamma_B$ such that
$\discr s$ preserves~$\CK_B$. A symmetry~$s$ is called
\emph{stable} if $\discr s$ acts identically on
$\CK_B^\perp/\CK_B$.

The group of symmetries (stable symmetries) of
the configuration $(\Gamma_B,\CK_B)$
is denoted by $\Sym(\Gamma_B,\CK_B)$ (respectively,
$\StSym(\Gamma_B,\CK_B)$\,).
\enddefinition

\Remark
In~\cite{JAG}, the configuration of a sextic~$B$ is defined
as the finite index
extension $\tS_B\supset S_B=\Sigma_B\oplus\<h_B\>$. In view of
Proposition~\ref{irreducible} and Theorem~\ref{th.Nik1}, in the
case of irreducible sextics the two definitions are equivalent.
\endRemark

\theorem\label{th.Dynkin}
For an irreducible plane sextic~$B$ with simple singularities,
the map $c\mapsto\tc_\Gamma$ establishes an isomorphism
$\StSym B\to\StSym(\Gamma_B,\CK_B)$.
\endtheorem

\proof
In view of Proposition~\ref{symplectic=>stable}, the exact sequence
given by Theorem~\ref{th.Sym} restricts to
$$
1@>>>\{\id,\tGt_*\}@>>>\Aut_\pm\CH_B@>>>\StSym B.
$$
Theorem~\ref{stable=>symplectic} provides a splitting
$c\mapsto\tc_*$ and, hence, an isomorphism
$$
\StSym B@>\cong>>\Aut_+\CH_B.
$$
Any element
$\tc_*\in\Aut_+\CH_B$ is uniquely determined by its
restriction to~$\Gamma_B$, and a symmetry
$\tc_\Gamma\in\Sym\Gamma_B$ extends to an element of
$\Aut_+\CH_B$ if and only if it is a stable symmetry of
$(\Gamma_B,\CK_B)$, see Corollary~\ref{extension}.
\endproof

\corollary\label{locally.constant}
Up to isomorphism,
the group $\StSym B$ depends only on the configuration
of~$B$. Furthermore, any path $B_t$, $t\in[0,1]$, in the
equisingular stratum $\CC_6(\Sigma_B)$ induces an isomorphism
$\StSym B_0\to\StSym B_1$.
\qed
\endcorollary

\section{Proof of Theorem~\ref{th.list}}\label{S.proof}

Throughout this section, $B$ is an irreducible plane sextic with
simple singularities. We use the notation introduced
in~\ref{s.K3}, abbreviating $\tX_B=\tX$, $\Gamma_B=\Gamma$, \etc.

\subsection{Sextics with type $\bE_8$ singular points}
Here, we treat the exceptional, in a certain sense, case
of curves that admit a stable symmetry but are not
\term2p-sextics.

\proposition\label{2E8}
Let $s\in\StSym(\Gamma_B,\CK_B)$ and $s\ne\id$. Then
either
\roster
\item
$B$ has two type~$\bE_8$ singular points, and $s$ is the
transposition of the two type~$\bE_8$ components of~$\Gamma_B$,
or
\item
$B$ is a \term2p-sextic, $p=3$, $5$, $7$, and $\discr s\ne\id$,
\endroster
the two cases being mutually exclusive.
\endproposition

\proof
Assume that $\discr s\ne\id$. Then, in order to make the action on
$\CK_B^\perp/\CK_B$ trivial, one must have $\CK_B\ne0$. According
to~\cite{degt.Oka}, $B$ is a \term2p-sextic, $p=3$, $5$, $7$; in
particular, $B$ has no type~$\bE_8$ singular points.

Now, assume that $\discr s=\id$. Then, in view of
Lemma~\ref{discr.Dynkin}, $s$ can only permute two
or more type~$\bE_8$ components of~$\Gamma_B$; in particular, $B$
has at least two type~$\bE_8$ singular points. On the other hand,
since $\mu(B)\le19$, the number of type~$\bE_8$ points is at most
two.
\endproof

Sextics with at least
two type~$\bE_8$ singular points are easily
classified using~\cite{JAG}; we merely state the final result.

\proposition\label{2E8.classes}
A sextic with
two type~$\bE_8$ singular points
can have one of the
five sets of singularities listed in Theorem~\iref{th.list}{E8}.
Each set of singularities is realized by a single
equisingular deformation family.
\qed
\endproposition

\subsection{\term2p-sextics}
Let $B$ be a \term2p-sextic, $p=3$, $5$, $7$. Then, according
to~\cite{degt.Oka}, the group $\CK_B\ne0$ is an $\F_p$-vector
space.
A singular
point~$P_i$ of~$B$ and the corresponding connected
component~$\Gamma_i$
of~$\Gamma_B$ is called \emph{essential} (\emph{ordinary}) if the
projection of~$\CK_B$ to $\discr\Gamma_i$ is non-zero
(respectively, zero).
Let $\bar\Gamma_B$ be the union of all essential components
of~$\Gamma_B$; it is obviously preserved by stable symmetries.
Denote by
$$
\pi_0\:\Sym\Gamma_B\to\SGSet(\pi_0(\Gamma_B))
\quad\text{and}\quad
\bar\pi_0\:\StSym(\Gamma_B,\CK_B)\to\SGSet(\pi_0(\bar\Gamma_B))
\eqtag\label{eq.pi}
$$
the corresponding representations in the symmetric groups.
Consider also the representation
$$
\Gk\:\StSym(\Gamma_B,\CK_B)\to\GL(\CK_B)\eqtag\label{eq.kappa}
$$
sending a symmetry
$s\in\StSym(\Gamma_B,\CK_B)$ to the restriction of $\discr s$ to~$\CK_B$.
Note that $\bar\pi_0(s)$ and $\Gk(s)$ are well defined on
$\Sym(\Gamma_B,\CK_B)$ and, hence, on
$\Aut_{\tGt}\tX_B$.
Furthermore, there is a homomorphism
$$
\bGk\:\Sym B\to\PGL(\CK_B)=\GL(\CK_B)/\pm\id\eqtag\label{eq.PGL}
$$
defined due to the fact that $\Gt_*$ induces $-\id$ on the
discriminant.


\proposition\label{fixed}
Let $B$ be a \term2p-sextic, $p=3$, $5$, $7$, and let
$s\in\StSym(\Gamma_B,\CK_B)$ be a stable symmetry of~$\Gamma_B$. Then\rom:
\roster
\item\local{ord}
$s$ acts trivially on $\Gamma_B\sminus\bar\Gamma_B$\rom;
\item\local{p}
$\discr s$ acts trivially on $\discr_q\Gamma_B$ for any prime
$q\ne p$\rom;
\item\local{fixed}
$s$ acts trivially on each component of~$\Gamma_B$ other than
$\bA_{p^r-1}$, $\bA_{2p^r-1}$, $r\ge1$, or~$\bE_6$
\rom(in the case $p=3$\rom)\rom;
\item\local{set}
$\pi_0(s)$ preserves each component of~$\Gamma_B$ other than
$\bA_{p^r-1}$, $r\ge1$, or~$\bE_6$ \rom(in the case $p=3$\rom).
\endroster
\endproposition

\proof
The first statement follows from the
fact that the discriminant of the union of the ordinary components
of~$\Gamma_B$ survives as a direct summand in $\CK_B^\perp/\CK_B$,
the fact that a \term2p-sextic has no singular points of
type~$\bE_8$, see~\cite{degt.Oka}, and Lemma~\ref{discr.Dynkin}.

Similarly, since $\CK_B\subset\discr_p\Gamma_B$, any other primary
component $\discr_q\Gamma_B$, $q\ne p$, survives to
$\CK_B^\perp/\CK_B$; this observation implies~\loccit{p}, and
the last two statements follow from~\loccit{p} and~\eqref{eq.discr}.
\endproof

\lemma\label{G.id}
Let $G$ be a finite abelian group and $s\:G\to G$ an automorphism
of order prime to $\ls|G|$. Assume that $G$ has an invariant
subgroup~$H$ such that the induced actions on~$H$ and $G/H$
are both trivial. Then $s=\id$.
\endlemma

\proof
For each element $g\in G$ one has $sg-g\in H$. Then, for some~$r$
prime to $\ls|G|$, one has
$g=s^rg=g+r(sg-g)$; hence, $r(sg-g)=0$ and $sg-g=0$.
\endproof

\corollary\label{CK.id}
Let $B$ be a \term2p-sextic, $p=3$, $5$, $7$.
Then the kernel of
the homomorphism $\Gk\:\StSym(\Gamma_B,\CK_B)\to\GL(\CK_B)$,
see \eqref{eq.kappa}, is a $p$-group.
\endcorollary

\proof
Let $s\in\StSym(\Gamma_B,\CK_B)$ be an element of order prime to~$p$, and
assume that $(\discr s)|_{\CK}=\id$.
Consider the filtration
$\discr_p\Sigma_B\supset\CK_B^\perp\supset\CK_B$. (Here, $^\perp$
stands for the orthogonal complement in $\discr_p\Sigma_B$.)
The action of $\discr s$ on~$\CK_B$ is trivial by the assumption,
the action on $\CK_B^\perp/\CK_B$ is trivial since
$s\in\StSym(\Gamma_B,\CK_B)$, and the action on
$\discr_p\Sigma_B/\CK_B^\perp$ is trivial due to the isomorphism
$\discr_p\Sigma_B/\CK_B^\perp=\Hom(\CK_B,\Q/\Z)$ given by the
discriminant bilinear form.
Applying Lemma~\ref{G.id} twice, we conclude that the action
of $\discr s$ on $\discr_p\Sigma_B$ is trivial. Hence, $s=\id$
due to Propositions~\iref{fixed}p and~\ref{2E8}.
\endproof

\corollary\label{F*}
Let $B$ be a \term2p-sextic, $p=3$, $5$, $7$, with
$\ell(\CK_B)=1$. Then the kernel of the homomorphism
$\Gk\:\StSym(\Gamma_B,\CK_B)\to\F_p^*=\GL(\CK_B)$
is a $p$-group.
\qed
\endcorollary

\corollary\label{fixed.component}
Let $B$ be
as in Corollary~\ref{F*},
and let $s\in\StSym(\Gamma_B,\CK_B)$ be an element of
order prime to~$p$. If the fixed point set of~$s$ contains an
essential component of~$\Gamma_B$, then $s=\id$.
\endcorollary

\proof
Under the assumption, $\Gk(s)=1\in\F_p^*$,
and Corollary~\ref{CK.id} applies.
\endproof

\corollary\label{Ker.pi}
If $B$ is a sextic as in Corollary~\ref{F*}, then
$\ls|\Ker\bar\pi_0|\le2$.
\endcorollary

\proof
In view of Proposition~\iref{fixed}{ord} and
Corollary~\ref{fixed.component}, any nontrivial element
$s\in\Ker\bar\pi_0$ has the following properties:
\roster
\item
the restriction of~$s$ to each ordinary component of~$\Gamma_B$
is~$\id$, and
\item
the restriction of~$s$ to each essential component of~$\Gamma_B$
is nontrivial.
\endroster
Since an essential component of~$\Gamma_B$ cannot be of
type~$\bD_4$ (\eg, due to~\eqref{eq.discr} and
Proposition~\ref{irreducible}),
it has at most one nontrivial symmetry, see
Lemma~\ref{Sym.Dynkin}; hence, the two properties above
determine~$s$ uniquely.
\endproof

\proposition\label{2.orbits}
Let $B$ be
as in Corollary~\ref{F*},
and let $s\in\StSym(\Gamma_B,\CK_B)$, $s\ne\id$,
be an element of order~$2$.
Then $\bar\pi_0(s)$ has at most two orbits.
\endproposition

\proof
Let~$V_p\subset\discr\Gamma_B$ be the $\F_p$-vector space of
order~$p$ elements, and let~$s_*$ be the
action of $\discr s$ on~$V_p$. Denote by~$V_p^-$
the $(-1)$ eigenspace of~$s_*$. Then, each orbit of $\bar\pi_0(s)$
contributes one to $\dim V_p^-$: each two element orbit
contributes a regular $\F_p$-representation of~$\CG2$, and each one
element orbit contributes a one dimensional representation $-\id$ due to
Corollary~\ref{fixed.component} and Lemma~\ref{-id.Dynkin}.
On the other hand, since $\dim\CK_B=1$, the
stability condition requires $\dim V_p^-\le2$.
\endproof

\corollary\label{Ker.pi=1}
Let $B$ be as in Corollary~\ref{F*}, and let
$\ls|\pi_0(\bar\Gamma_B)|\ge3$. Then $\Ker\bar\pi_0=1$.
\qed
\endcorollary

\subsection{\term10-sextics}\label{s.D10}
All \term10-sextics are classified in~\cite{degt.Oka}: they form
eight equisingular deformation families, their sets of essential
singularities are $4\bA_4$, $\bA_9\splus2\bA_4$, or $2\bA_9$, and
for any such sextic~$B$ one has $\CK_B=\Z_5$.

Any symmetry of~$\Gamma_B$ of order divisible by~$5$ would have to
permute cyclically at least five isomorphic components
of~$\Gamma_B$. Hence, such a symmetry does not exist, and
Corollary~\ref{F*} implies that
$\StSym(\Gamma_B,\CK_B)\subset\F_5^*\cong\CG4$.
Symmetries of order~$2$ are
constructed in~\cite{degt.Oka3}, and to show that
$\StSym(\Gamma_B,\CK_B)\cong\CG2$,
it remains to rule out symmetries of order~$4$.

Let $s\in\StSym(\Gamma_B,\CK_B)$ be an element of order~$4$. Applying
Corollary~\ref{fixed.component} to~$s^2$
and using Lemma~\ref{Sym.Dynkin}, one concludes that
$\bar\pi_0(s)$ has no fixed points.
Then,
due to Proposition~\iref{fixed}{set}, the essential
singularities of~$B$ are $4\bA_4$, and Proposition \ref{2.orbits}
applied to~$s^2$ implies that $\bar\pi_0(s)$ is a cycle of length~$4$.
Thus, $\discr s$ is a
regular $\F_5$-representation of~$\CG4$; it has four distinct
eigenspaces (one for each eigenvalue $\Gl\in\F_5^*$)
and thus cannot be stable.

\subsection{\term14-sextics}\label{s.D14}
According to~\cite{degt.Oka}, \term14-sextics form two
equisingular deformation families; the set of essential
singularities of any such sextic~$B$ is $3\bA_6$, and one has
$\CK_B=\Z_7$.

As in~\ref{s.D10}, the graph~$\Gamma_B$ has no symmetries of
order~$7$ and Corollary~\ref{F*} implies that
$\StSym(\Gamma_B,\CK_B)\subset\F_7^*\cong\CG6$.
Symmetries of order~$3$ are
constructed in~\cite{degt-Oka}. Assume that
$s\in\StSym(\Gamma_B,\CK_B)$ is
a nontrivial element of order~$2$.
Then, due to Proposition~\ref{2.orbits}
and Corollary~\ref{fixed.component}, $\discr s$ is the direct sum
of a dimension one representation $-\id$
and a regular $\F_7$-representation of~$\CG2$; it is
easy to see that this action has no isotropic invariant
subspaces.

\subsection{\term6-sextics: symmetries of order~$3$}
The sets of singularities of \term6-sextics (\ie,
irreducible sextics of torus type) are classified in
M.~Oka, D.~T.~Pho~\cite{OkaPho.moduli}. If $w(B)\le7$, the set of
essential singularities of~$B$ is of the form
$$
\tsize
\bigsplus_ik_i\bA_{3i-1}\splus l\bE_6,\quad
\sum_iik_i+2l=6,\eqtag\label{eq.w=6}
$$
and one has $\CK_B=\Z_3$. (In this case, essential
are the inner singular points with respect to the only torus
structure of~$B$.)
If $w(B)=8$, then the set of essential singularities
is
$$
\bE_6\splus\bA_5\splus4\bA_2,\quad
\bE_6\splus6\bA_2,\quad
2\bA_5\splus4\bA_2,\quad
\bA_5\splus6\bA_2,\quad\text{or}\quad
8\bA_2,\eqtag\label{eq.w=8}
$$
and one
has $\CK_B=\Z_3\oplus\Z_3$. Finally, there is one deformation
family of sextics of weight~$9$: their set of singularities is
$9\bA_2$, all nine cusps being essential,
and one has $\CK_B=\Z_3\splus\Z_3\splus\Z_3$. Stable
symmetries of sextics of weight~$9$ are described
in~\cite{degt.8a2}; for this reason, we ignore them here.

Let $B$ be a plane sextic with simple singularities (not
necessarily a \term6-sextic or even irreducible),
and let $\tc\in\Aut_{\tGt}\tX_B$
be a symplectic automorphism of order~$3$.
Denote by~$\tX'$
the minimal resolution of singularities of the quotient
$\tX_B/\tc$; it is also a $K3$-surface (since $\tc$ is
symplectic). Let $\Gamma'$ be the union of the
components of~$\Gamma_B$ fixed by~$\pi_0(\tc_\Gamma)$.

\lemma\label{lem.order3}
In the notation above, $\Gamma'$ has the form
$$\tsize
\bigocup_ik_i\bA_i\ocup l\bD_4,\quad\text{where}\quad
f=\sum_ik_i(i+1)+2l\le6.
$$
The image in~$\tX'$ of the divisor represented by each
type~$\bA_i$ \rom(respectively,~$\bD_4$\rom)
component of~$\Gamma'$ is a union
of $(-2)$-curves in~$\tX'$ whose incidence graph is $\bA_{3i+2}$
\rom(respectively,~$\bE_6$\rom). In addition, $\tc$ has $(6-f)$
isolated fixed points not in the $(-2)$-curves represented by the
vertices of~$\Gamma_B$.
\endlemma

\proof
%
Let $E_j\subset\tX_B$ be the divisor represented by a component
$\Gamma_j$ fixed by~$\tc_\Gamma$ pointwise.
$E_j$ is a union of $(-2)$-curves, and each
point of intersection of two distinct
components of~$E_j$ is a fixed point
of~$\tc$. On the other hand, any $(-2)$-curve preserved
by~$\tc$
contains exactly two fixed points of~$\tc$. (Since $\tc$ is
symplectic, its fixed points are isolated.)
Hence, $\Gamma_j$
cannot be of type~$\bD$ or~$\bE$ (as otherwise $E_j$ would have a
component with three fixed points, hence fixed pointwise),
and if it is of type~$\bA_i$,
then $E_j$ contains $(i+1)$ fixed points of~$\tc$.

Now,
let $\Gamma_j\subset\Gamma'$ be a component \emph{not} fixed
by~$\tc_\Gamma$ pointwise.
According to Lemma~\ref{Sym.Dynkin}, it must be of type~$\bD_4$,
and the divisor~$E_j$ represented by~$\Gamma_j$
has a single component (corresponding
to the central vertex of~$\Gamma_j$) fixed by~$\tc$; this
component contains two fixed
points of~$\tc$.

As is well known (see, \eg,
Nikulin~\cite{Nikulin.auto}),
any symplectic automorphism of order~$3$ has six fixed points.
Combining this observation with the count above, one obtains the
estimate on the number and types of the components of~$\Gamma'$.
The calculation
of the image in~$\tX'$ is straightforward: each fixed point
of~$\tc$ gives rise to a cusp in $\tX_B/\tc$, \ie, to two
extra $(-2)$-curves in~$\tX'$.
\endproof

\corollary\label{cor.order3}
Let~$B$ be a \term6-sextic with $w(B)\le8$,
and let $\tc\in\Aut_{\tGt}\tX_B$ be a symplectic automorphism
of order~$3$.
Then the set of
orbits of $\bar\pi_0(\tc_\Gamma)$ is one of the following\rom:
$(3\bE_6)$, $(3\bA_5)$, or $(3\bA_2)\ocup(3\bA_2)$.
\endcorollary

\proof
If $w(B)\le7$, the statement follows immediately
from Lemma~\ref{lem.order3} and \eqref{eq.w=6}.
Assume that $w(B)=8$. Due to Lemma~\ref{lem.order3}
and~\eqref{eq.w=8}, the orbits of $\bar\pi_0(\tc_\Gamma)$ are
either $(\bA_2)\ocup(\bA_2)\ocup(3\bA_2)\ocup(3\bA_2)$ or
$(\bA_5)\ocup(3\bA_2)\ocup(3\bA_2)$.
Hence, the projection of the exceptional divisor to~$\tX'$ is
a divisor $E\subset\tX'$ with the incidence
graph $2\bA_8\ocup\bA_2\ocup\bA_2$ or
$\bA_{17}\ocup\bA_2\ocup\bA_2$.
It spans a negative definite lattice of rank $\ge20$,
which does not fit into $H_2(\tX')$.
\endproof

\corollary\label{3E6.order3}
A \term6-sextic of weight~$\le8$
has a stable symmetry of order~$3$ if and only if
its set of essential singularities is $3\bE_6$.
\endcorollary

\proof
A cyclic permutation of three type~$\bA_5$ components cannot be
stable due to Proposition~\iref{fixed}{set}. For two cycles of
length~$3$ on six cusps, the induced action on $\discr\Gamma_B$ is
a direct sum of two copies of a regular $\F_3$-representation
of~$\CG3$; it cannot be stable since $\dim\CK_B=1$. The existence
of a stable symmetry of order~$3$ on the set of essential
singularities~$3\bE_6$ is obvious.
\endproof

\subsection{\term6-sextics of weight $\le7$}
Propositions~\ref{w=6.order2}
below describes the image of the
representation $\Gk\:\StSym(\Gamma_B,\CK_B)\to\F_3^*=\Z_2$, see
Corollary~\ref{F*}.

\proposition\label{w=6.order2}
Let $B$ be a \term6-sextic, $w(B)\le7$, and let
$s\in\StSym(\Gamma_B,\CK_B)$
be a nontrivial element of order~$2$. Then the
orbits of~$\bar\pi_0(s)$ are as follows\rom:
$$
(2\bE_6)\ocup(\bE_6),\quad
(2\bE_6)\ocup(\bA_5),\quad
(2\bE_6)\ocup(2\bA_2),\quad
(\bA_{17}),\quad
(2\bA_8).
$$
Conversely, any set of orbits as above is realized by a
stable symmetry of order~$2$.
\endproposition

\proof
In view of Proposition~\iref{fixed}{fixed},~\ditto{set},
each orbit is one of the following: $2\bE_6$, $\bE_6$, $\bA_{17}$,
$2\bA_8$, $\bA_8$, $\bA_5$, $2\bA_2$, or~$\bA_2$. Due to
Proposition~\ref{2.orbits}, there are at most two orbits.
Combining these
observations with~\eqref{eq.w=6}, one obtains the five sets of
orbits
listed in the statement and $(\bA_8)\ocup(\bA_8)$. In the latter
case, disregarding the ordinary components, one has
$\discr s=-\id$, see Corollary~\ref{fixed.component}
and Lemma~\ref{-id.Dynkin}, and since
$\CK_B^\perp/\CK_B\cong\Z_9$, the symmetry is not stable.

The converse statement is straightforward:
all five involutions are easily constructed
using the description of~$\CK_B$
given in~\cite{degt.Oka}.
\endproof

\corollary\label{w=6.kappa}
Let $B$ be a \term6-sextic of weight $w(B)\le7$. The
representation $\Gk\:\StSym B\to\F_3^*=\Z_2$ is an
epimorphism
if and only if the set~$\Sigma_B^{\inj}$ of
essential singularities of~$B$ is as in
Theorem~\iref{th.list}{3E}, \ditto{E}, or~\ditto{A}.
The kernel $\Ker\kappa$ is trivial unless
$\Sigma_B^{\inj}=3\bE_6$, \ie, unless
$B$ is as in Theorem~\iref{th.list}{3E}.
\endcorollary

\proof
The epimorphism part follows from Proposition~\ref{w=6.order2},
and the kernel of~$\kappa$ is estimated using
Corollaries~\ref{F*} and~\ref{3E6.order3}.
\endproof

\proposition\label{3E6}
For a \term6-sextic~$B$ with
$\Sigma_B^{\inj}=3\bE_6$, \cf\. Theorem~\iref{th.list}{3E},
the representation
$\bar\pi_0\:\StSym(\Gamma_B,\CK_B)\to\SGSet(\bar\Gamma_B)=\SG3$ is an
isomorphism.
\endproposition

\proof
$\bar\pi_0$ is an epimorphism due to
Proposition~\ref{w=6.order2} and Corollary~\ref{3E6.order3}.
Its kernel is trivial due to Corollary~\ref{Ker.pi=1}.
\endproof

\subsection{Sextics of weight eight}\label{s.w=8}
Let $B$ be a \term6-sextic of weight~$8$. According to
Corollaries~\ref{CK.id} and~\ref{cor.order3}, the representation
$\Gk\:\StSym(\Gamma_B,\CK_B)\to\GL(\CK_B)$
as in~\eqref{eq.kappa}
is monic. In~\cite{degt.8a2}, we constructed a stable
symmetry~$c$ of order~$2$ whose image in $\GL(\CK_B)$ is the
central element $-\id$; the minimal resolution of singularities of
the quotient $\Cp2\!/c$ is a geometrically ruled rational
surface~$\Sigma_2$ with an exceptional section~$E$ of
self-intersection~$(-2)$, and the image $B/c$ is a trigonal curve
$\B\subset\Sigma_2$ with four cusps,
\cf. Section~\ref{s.involutions} and Theorem~\ref{th.maximal}.
The original plane~$\Cp2$ is
the double covering of the quadratic cone $\Sigma_2/E$ ramified
at the vertex $E/E$ and a certain section~$\L$ of~$\Sigma_2$, and
$B$ is the pull-back of~$\B$.

Since $c$ is central, any other stable symmetry of~$B$ would
descend to a symmetry of $(\Sigma_2,\B+\L)$ stable under
equisingular deformations of $\B+\L$.
The curve~$\B$ is rigid;
in appropriate affine coordinates $(x,y)$ in~$\Sigma_2$ it is
given by the polynomial
$$
f(x,y)=4y^3-(24x^3+3)y+(8x^6+20x^3-1).
$$
The group of (Klein) symmetries of~$B$ is the alternating
group~$\AG4$ (respectively, symmetric group~$\SG4$); it can
be identified with the group of even (respectively, all)
permutations of the four
cusps of~$\B$, see~\cite{degt.8a2}. (Recall that a \emph{Klein
automorphism} of an analytic variety is an either holomorphic or
anti-holomorphic automorphism.)
Since $B$ has no stable
symmetries of order~$3$, it suffices to show that
$\B+\L$ has no stable symmetries of order~$2$. (Note that
$\ls|\GL(\CK_B)|=(3)(2)^4$.)
All order~$2$
elements in~$\AG4$ are conjugate, and one of them is given by the
change of coordinates
$$
(x,y)=\Bigl(-\frac{x'-\Ge}{2\Ge^2x'+1},-\frac{3y'}{(2\Ge^2x'+1)^2}\Bigr),
\quad\Ge=-\frac12+i\frac{\sqrt3}2.
$$
A section of~$\Sigma_2$ is preserved by this transformation if and
only if it has the form $y=a(2x^2+2\Ge x-\Ge^2)$, $a\in\C$,
and it is
straightforward that the family $\B+\L$ with~$\L$ as above does
not contain any equisingular stratum; hence, $B$ has no other
stable symmetries.

\Remark
There is a unique section $\L=\{y=0\}$ that is invariant under the
full group~$\SG4$ of Klein
symmetries of~$\B$. It gives rise to a unique, up to projective
transformation, sextic~$B$ of weight~$8$ (with the set of
singularities $8\bA_2$) for which the image of the representation
$\bGk\:\Sym B\to\PGL(\CK_B)\cong\SG4$, see~\eqref{eq.PGL},
is the subgroup $\AG4\subset\SG4$. Using~\cite{degt.8a2}, this
image can be identified with the group of even permutations of the
four torus structures of~$B$.
\endRemark

\subsection{Proof of Theorem~\ref{th.list}}\label{proof.list}
Proposition~\ref{2E8} states that most irreducible sextics
admitting stable symmetries are \term2p-sextics, $p=3$, $5$,~$7$.
The exceptional case of
sextics with two type~$\bE_8$ singular points is covered by
Proposition~\ref{2E8.classes}.

The group of stable symmetries of \term10- and \term14-sextics are
found in Sections~\ref{s.D10} and~\ref{s.D14}, respectively.
\term6-sextics of weight $w\le7$ are covered by
Corollary~\ref{w=6.kappa} and Proposition~\ref{3E6}, sextics of
weight~$8$ are considered in Section~\ref{s.w=8}, and the
remaining case of sextics of weight~$9$ is contained
in~\cite{degt.8a2}.
\qed

\section{Stable involutions}\label{S.involutions}

In this section, we analyze the relation between stable
involutions of irreducible sextics and maximal trigonal curves
in~$\Sigma_2$.

\subsection{Maximal trigonal curves in $\Sigma_2$}\label{s.maximal}
The \emph{Hirzebruch surface} $\Sigma_k$, $k\ge0$,
is a geometrically ruled
rational surface with an exceptional section~$E$ of square~$-k$. A
\emph{trigonal curve} is a curve $B\subset\Sigma_k$ disjoint
from~$E$ and
intersecting each generic fiber at three points. In appropriate
affine coordinates $(x,y)$ in~$\Sigma_k$, a trigonal curve can be
given by its \emph{Weierstra{\ss} equation}
$y^3+g_2(x)y+g_3(x)=0$,
where $\deg g_2\le2k$ and $\deg g_3\le3k$, and the
\emph{\rom(functional\rom) $j$-invariant} of~$B$ is defined as the
function
$$
j=j_B\:\Cp1\to\Cp1,\quad
x\mapsto\frac{4g_2^3(x)}{\Delta(x)},\quad
\text{where $\Delta=4g_2^3+27g_3^2$}.\eqtag\label{eq.j}
$$
Here, the first copy of~$\Cp1$ (the source) is the base of the
ruling of~$\Sigma_k$, and the second copy (the target) is the
standard Riemann sphere $\C\cup\infty$. The curve~$B$ is called
\emph{isotrivial} if $j_B=\const$.

By a \emph{singular fiber} of a trigonal curve $B\subset\Sigma_k$
we mean a fiber of the ruling of~$\Sigma_k$ intersecting~$B$
geometrically at less than three points. Locally, in a
neighborhood of a simple singular fiber, $B$ is the
ramification locus of the Weierstra\ss{} model of a Jacobian
elliptic surface, and to describe the type of the fiber we use
(one of) the standard notation for the singular elliptic fibers,
referring to the extended Dynkin graph of the exceptional
divisors. For non-simple singular fibers, we use Arnol$'$d's
notation $\tilde\bJ_{k,p}$ and $\tilde\bE_{6k+\epsilon}$, $k\ge2$,
referring to the type of the singular point of~$B$.

The $j$-invariant has three special values, $0$, $1$,
and~$\infty$, which are typically taken at the roots of~$g_2$,
$g_3$, and~$\Delta$, respectively. In~\cite{degt.kplets}, a
trigonal curve~$B$ with double singular points only
is called \emph{maximal} if
$j_B$ has the following properties:
\roster
\item\local1
$j_B$ has no critical values other than~$0$, $1$, or~$\infty$, and
\item\local2
each pull-back $j_B^{-1}(0)$ (respectively,~$j_B^{-1}(1)$)
has ramification index at most three (respectively, at most two).
\endroster
In order to extend this definition to all trigonal curves, we need
to exclude singular fibers similar to~$\tilde\bD_4$, which
are not detected by the $j$-invariant and typically increase the
dimension of the moduli space. (Essentially, the additional
requirement means that each singular fiber should remain singular
after elementary transformations.)
Thus, we have the following definition.

\definition\label{def.maximal}
A trigonal curve $B\subset\Sigma_k$ is called \emph{maximal} if
$B$ has no singular fibers of type~$\tilde\bD_4$ or
$\tilde\bJ_{k,0}$, $k\ge2$, and
$j_B$ satisfies
conditions~\itemref{s.maximal}1,~\ditto2 above.
\enddefinition

With this definition, the alternative characterization of maximal
curves given in~\cite{degt.kplets} still holds: {\proclaimfont
a non-isotrivial curve~$B$
is maximal if and only if it does not admit a nontrivial
degeneration to another non-isotrivial trigonal curve.}

A singular fiber of~$B$ is called \emph{stable} if it is preserved
by small equisingular (but not necessarily fiberwise) deformations
of~$B$. Clearly, stable are all fibers except those of
type~$\tilde\bA_0^{**}$, $\tilde\bA_1^*$, or~$\tilde\bA_2^*$
(which can split into a stable fiber and $\tilde\bA_0^*$).
The curve~$B$ is
called \emph{stable} if all its singular fibers are stable.

\theorem\label{th.stable.curves}
Up to automorphism of~$\Sigma_2$, a stable maximal trigonal curve
$B\subset\Sigma_2$ is determined by its set of singular fibers,
which can be one of
those listed in Table~\ref{tab.curves}.
\midinsert
\table\label{tab.curves}
Singular fibers of stable maximal curves in~$\Sigma_2$
\endtable
\def\col#1#2{\vtop\bgroup\halign\bgroup$##$\hfil&&\quad\ $##$\hfil\cr
 \multispan#1\hidewidth\rm#2\hidewidth\cr\noalign{\smallskip}}
\def\endcol{\crcr\egroup\egroup}
\vskip-2pt
\centerline{\openup2pt
\col2{Irreducible curves}
\tilde\bE_8\splus2\tilde\bA_0^*&
 2\tilde\bA_4\splus2\tilde\bA_0^*\cr
\tilde\bE_6\splus\tilde\bA_2\splus\tilde\bA_0^*&
 4\tilde\bA_2\cr
\tilde\bA_8\splus3\tilde\bA_0^*\cr
\endcol\qquad\quad
\col2{Reducible curves}
\tilde\bE_7\splus\tilde\bA_1\splus\tilde\bA_0^*&
 \tilde\bA_7\splus\tilde\bA_1\splus2\tilde\bA_0^*\cr
\tilde\bD_8\splus2\tilde\bA_0^*&
 \tilde\bA_5\splus\tilde\bA_2\splus\tilde\bA_1\splus\tilde\bA_0^*\cr
\tilde\bD_6\splus2\tilde\bA_1&
 2\tilde\bA_3\splus2\tilde\bA_1\cr
\tilde\bD_5\splus\tilde\bA_3\splus\tilde\bA_0^*\cr
\endcol}
\endinsert
\endtheorem

\Remark\label{rem.isotrivial}
In Table~\ref{tab.curves}, the curves with a type~$\tilde\bE$
singular fiber (and only these curves) admit equisingular isotrivial
degenerations: $\tilde\bE_8\splus\tilde\bA_0^{**}$,
$\tilde\bE_7\splus\tilde\bA_1^{*}$, and
$\tilde\bE_6\splus\tilde\bA_2^{*}$.
\endRemark

\proof
According to~\cite{degt.kplets}, a maximal trigonal curve
$B\subset\Sigma_k$ with
double singular points only
is determined, up to automorphism
of~$\Sigma_k$, by its \emph{skeleton} $\Sk B\subset\Cp1\cong S^2$
(Grothendieck's \emph{dessin d'enfants}), which is defined as
the bi-partite planar map $j_B^{-1}([0,1])$, the pull-backs of~$0$
and~$1$ being, respectively, black and white vertices of~$\Sk B$.
The skeleton has the following properties:
\roster
\item\local1
$\Sk B$ is connected;
\item\local2
$\Sk B$ has at least one black and at least one white vertex;
\item\local3
the valency of each black (white) vertex is~$\le3$
(respectively,~$\le2$).
\endroster
Conversely, ani bi-partite planar map
satisfying~\loccit1--\loccit3
above is the skeleton of a certain maximal trigonal curve
$B\subset\Sigma_k$; the parameter~$k$ is given by the relation
$2k=b_1+2b_2+b_3+w_1$, where $b_i$ and $w_i$ are the numbers of,
respectively, black and white vertices of valency~$i$.

Since the skeleton~$\Sk B$ is a bi-partite graph, each
complementary region of $\Sk B$ has equal numbers of black and
white vertices (`corners') in the boundary; we call a region a
\emph{$p$-gon}, $p\ge1$, if it has $p$ black and $p$ white
corners. The stable singular fibers of~$B$ are in a one-to-one
correspondence with the regions of~$\Sk B$:
each $p$-gonal region
contains a single singular fiber of type $\tilde\bA_{p-1}$
($\tilde\bA_0^*$ if $p=1$). The unstable fibers of
type~$\tilde\bA_0^{**}$, $\tilde\bA_1^*$, and~$\tilde\bA_2^*$ are
over, respectively,
the $1$-valent black vertices, $1$-valent white
vertices, and $2$-valent black vertices.
For this reason, we call black vertices of valency $\le2$ and
white vertices of valency~$1$ \emph{unstable}.

\midinsert
\centerline{\vbox{\halign{\hss#\hss&&\qquad\hss#\hss\cr
\cpic{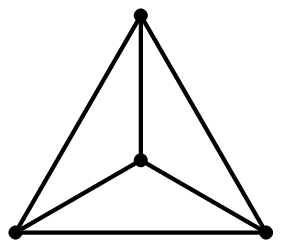}&\cpic{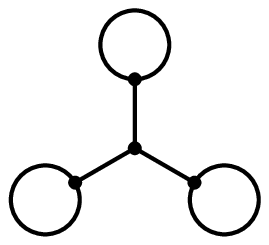}&\cpic{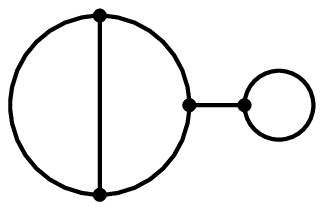}\cr
\noalign{\medskip}
$4\tilde\bA_2$&$\tilde\bA_8\splus3\tilde\bA_0^*$&
 $\tilde\bA_5\splus\tilde\bA_2\splus\tilde\bA_1\splus\tilde\bA_0^*$\cr
\crcr}}}
\bigskip
\centerline{\vbox{\halign{\hss#\hss&&\qquad\hss#\hss\cr
\cpic{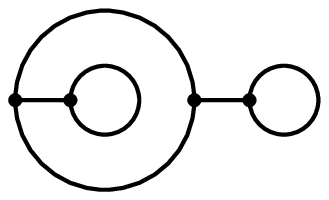}&\cpic{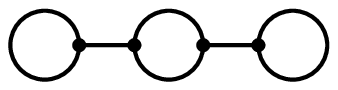}&\cpic{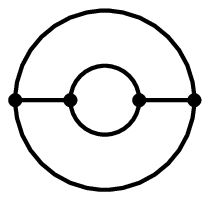}\cr
\noalign{\medskip}
$2\tilde\bA_4\splus2\tilde\bA_0^*$&
 $\tilde\bA_7\splus\tilde\bA_1\splus2\tilde\bA_0^*$&
 $2\tilde\bA_3\splus2\tilde\bA_1$\cr
\crcr}}}
\figure\label{fig.Sigma2}
Stable maximal curves in $\Sigma_2$ with double points only
\endfigure
\endinsert

Thus, in order to classify stable trigonal curves in~$\Sigma_2$
with double singular points only,
it suffices to list al bi-partite planar maps
satisfying~\loccit1--\loccit3 above, with four trivalent black
vertices, and without unstable vertices. This is done in
Figure~\ref{fig.Sigma2}. (In the figures, we omit bivalent white
vertices; one such vertex is to be placed at the center of each
edge connecting two black vertices.) Reducible curves are detected
using the criterion found in~\cite{degt.kplets}. (It is worth
mentioning that the skeleton is a graph in the \emph{oriented}
sphere~$\Cp1$. However, all graphs shown in
Figures~\ref{fig.Sigma2} and~\ref{fig.Sigma1} are symmetric. In
particular, this means that all curves are real.)

If the curve has one simple triple point, we apply an elementary
transformation centered at this point and obtain a maximal
trigonal curve $B_1\subset\Sigma_1$ with at most one unstable
fiber. Such curves are classified in Figure~\ref{fig.Sigma1}. If
$B_1$ is stable, the inverse elementary transformation can
contract any singular fiber of~$B$,
resulting in the curves with a type~$\tilde\bD$
singular fiber in Table~\ref{tab.curves}. If $B_1$ has an unstable
singular fiver of
type~$\tilde\bA_0^{**}$, $\tilde\bA_1^*$, or~$\tilde\bA_2^*$, the
inverse elementary transformation should contract this fiber,
resulting in a singular fiber of type~$\tilde\bE_6$,
$\tilde\bE_7$, or~$\tilde\bE_8$, respectively. To detect the
reducible curves, one can either use the criterion
in~\cite{degt.kplets} or just notice that trigonal curves
in~$\Sigma_1$ are merely plane cubics.

\midinsert
\centerline{\vbox{\halign{\hss#\hss&&\kern1.5em \hss#\hss\cr
\cpic{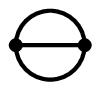}&\cpic{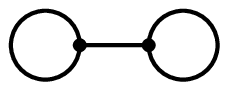}&\cpic{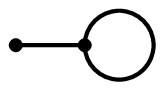}&\cpic{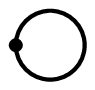}&\cpic{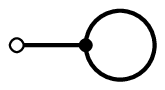}\cr
\noalign{\medskip}
$3\tilde\bA_1$&$\tilde\bA_3\splus3\tilde\bA_0^*$&
 $\tilde\bA_0^{**}\splus\tilde\bA_2\splus\tilde\bA_0^*$&
 $\tilde\bA_2^*\splus2\tilde\bA_0^*$&
 $\tilde\bA_1^*\splus\tilde\bA_1\splus\tilde\bA_0^*$\cr
\crcr}}}
\figure\label{fig.Sigma1}
Maximal curves in $\Sigma_1$ with at most one unstable fiber
\endfigure
\endinsert

To complete the proof, it remains to notice that a non-isotrivial
trigonal curve in~$\Sigma_2$ cannot have two triple points or a
non-simple triple point (adjacent to~$\bJ_{10}$ in Arnol$'$d's
notation), as otherwise one would apply two elementary
transformations and obtain a trigonal curve in~$\Sigma_0$, which
is necessarily isotrivial.
\endproof

\corollary\label{cor.stable.curves}
Up to automorphism of~$\Sigma_2$, a stable maximal trigonal curve
$B\subset\Sigma_2$ is determined by its set of singularities.
\endcorollary

\proof
The set of singularities of~$B$ is obtained from its set of
singular fibers by disregarding the type~$\tilde\bA_0^*$ summands
and `removing the tildes.' From Table~\ref{tab.curves} it follows
that the two sets determine each other. Furthermore, the
maximality of a set of singular fibers
can be tested numerically, by applying the
Riemann--Hurwitz formula to~$j_B$. (For example, if all singular
fibers are of type~$\tilde\bA_p$, $p\ge1$, or~$\tilde\bA_0^*$, the
curve is maximal if and only if the number of singular fibers is
four.) Hence, each set of singularities obtained from
Table~\ref{tab.curves} is realized by maximal curves
(or their equisingular isotrivial degenerations,
see Remark~\ref{rem.isotrivial}) only.
\endproof

\corollary\label{mu=8}
A non-isotrivial trigonal curve $B\subset\Sigma_8$ is stable and
maximal if and only if $\mu(B)=8$.
\endcorollary

\proof
The direct statement follows from Table~\ref{tab.curves}. For the
converse, we compare two independent classifications. A necessary
condition for a set of simple singularities~$\Sigma$ to be
realized by a trigonal curve in~$\Sigma_2$ is that $\Sigma$,
regarded as a root system, must admit an embedding to~$\bE_8$,
see, \eg,~\cite{degt.Oka2}. In addition to those listed in
Table~\ref{tab.curves}, there are three such root systems of rank
eight: $2\bD_4$, $\bD_4\splus4\bA_1$, and $8\bA_1$. The former set
of singularities is realized by isotrivial curves (obtained by
two elementary transformations from a union of three disjoint
sections of~$\Sigma_0$). The sets of singularities
$\bD_4\splus4\bA_1$ and $8\bA_1$ are not realized by a trigonal
curve in~$\Sigma_2$: each singular fiber of type~$\tilde\bD_4$
(respectively, $\tilde\bA_1$) is a root of
the discriminant~$\Delta$, see~\eqref{eq.j}, of
multiplicity~$6$ (respectively,~$2$), whereas $\deg\Delta\le12$.
\endproof

\Remark
As it follows from the proof of Corollary~\ref{mu=8},
the isotrivial curves
$B\subset\Sigma_2$ with $\mu(B)=8$ are either the equisingular
isotrivial
degenerations listed in Remark~\ref{rem.isotrivial} or the curve
with the set of singular fibers $2\tilde\bD_4$.
\endRemark

\subsection{Stable involutions}\label{s.involutions}
Recall that the set of fixed points of an involutive
automorphism~$c$ of~$\Cp2$ consists of a line~$L_c$ and an
isolated point~$O_c$. The blown up quotient $\Cp2(O_c)/c$ is the
Hirzebruch surface~$\Sigma_2$;
the exceptional divisor
over~$O_c$ projects to the exceptional section $E\subset\Sigma_2$,
and the line~$L_c$ projects to a
generic section $\L\subset\Sigma_2$ disjoint from~$E$.
Conversely, given a section $\L\subset\Sigma_2$ disjoint from~$E$,
the double covering of $\Sigma_2/E$ ramified at~$\L$ and $E/E$ is
the plane~$\Cp2$, and the deck translation of the covering is an
involutive automorphism whose fixed point set is the pull-back of
the union $\L\cup(E/E)$.

\theorem\label{th.maximal}
Let $B\subset\Cp2$ be an irreducible sextic with simple
singularities, and let $c\in\StSym B$
be an involutive stable symmetry of~$B$.
Then the image of~$B$ in the Hirzebruch surface
$\Sigma_2=\Cp2(O_c)/c$ is an irreducible stable maximal trigonal
curve~$\B$ \rom(or an equisingular
isotrivial degeneration of such a curve, see
Remark~\ref{rem.isotrivial}\rom).
The set of singularities of~$\B$ is as follows\rom:
\roster
\item
$\bE_8$, if $B$ is as in~\iref{th.list}{E8}\rom;
\item
$\bE_6\splus\bA_2$, if $B$ is as in~\iref{th.list}{3E}
or~\ditto{E}\rom;
\item
$\bA_8$, if $B$ is as in~\iref{th.list}{A}\rom;
\item
$2\bA_4$, if $B$ is as in~\iref{th.list}{10}\rom;
\item
$4\bA_2$, if $B$ is as in~\iref{th.list}{w=9}
or~\ditto{w=8}.
\endroster
\endtheorem

\proof
According to~\cite{degt.Oka3}, the image~$\B$ is either a trigonal
curve or a hyperelliptic curve with $\B\circ E=2$; in both cases,
the singularities of~$\B$ can be found using the results
of~\cite{degt.Oka3}. Assuming that $\B$ is a trigonal curve, the
essential singular points of~$B$ project to the sets of
singularities listed in the statement, while the ordinary singular
points give rise to points of tangency of~$\B$ and~$\L$. To
complete the proof in this case, one applies Corollary~\ref{mu=8}.

It remains to rule out the possibility that $\B$ is a
hyperelliptic curve. As, in this case,
$\B$ cannot have a triple point, it cannot appear from a
sextic~$B$ as in~\iref{th.list}{3E}, \ditto{E}, or~\ditto{E8}. The
sextics as
in~\iref{th.list}{w=9},~\ditto{w=8} and~\iref{th.list}{10} were
treated in~\cite{degt.8a2} and~\cite{degt.Oka3}, respectively. The
only remaining possibility is the set of essential singularities
$\bA_{17}$ in~\iref{th.list}{A}, which can be obtained from a
hyperelliptic curve~$\B$ with a single type~$\bA_7$ singular
point on~$E$. Such a curve~$\B$ does exist, but it is necessarily
reducible (see, \eg,~\cite{quintics}); hence, so is~$B$.
\endproof

\theorem\label{th.converse}
An involutive symmetry~$c$ of an irreducible plane sextic~$B$ with
simple singularities is stable if and only if the image of~$B$
in the Hirzebruch surface
$\Sigma_2=\Cp2(O_c)/c$ is an irreducible stable maximal trigonal
curve \rom(or an equisingular
isotrivial degeneration of such a curve, see
Remark~\ref{rem.isotrivial}\rom).
\endtheorem

\proof
The `only if' part is given by Theorem~\ref{th.maximal}. The `if'
part follows essentially from comparing Theorems~\ref{th.maximal}
and~\ref{th.stable.curves}: in addition, one needs to check that,
for each degeneration of the section~$\L$ (passing through a
singular point of~$\B$, tangency to~$\B$, \etc.),
the dimension of the
moduli space of such sections coincides with its expected
dimension, which in turn equals the dimension of the moduli space
of corresponding sextics, \cf. Remark~\ref{conjecture.proof}
below.
We leave details to the reader.
(In fact, the sets of singularities $2\bA_4$ and $4\bA_2$ were
studied in~\cite{degt.Oka3} and~\cite{degt.8a2}, respectively; the
three other sets of singularities will be considered in a
subsequent paper.)
\endproof

\conjecture\label{conjecture}
An involutive symmetry~$c$ of a plane sextic~$B$ with
simple singularities only is stable if and only if the image of~$B$
in the Hirzebruch surface
$\Sigma_2=\Cp2(O_c)/c$ is a stable maximal trigonal
curve \rom(or an equisingular
isotrivial degeneration of such a curve, see
Remark~\ref{rem.isotrivial}\rom).
\endconjecture

\Remark\label{non-simple}
A simple parameter count shows that, for sextics with a non-simple
singular point, the conclusion of Conjecture~\ref{conjecture}
fails. For example, some sextics~$B$ with the sets of singularities
$\bY_{1,1}^1\splus2\bA_4$, $\bY_{1,1}^1\splus\bA_9$, and
$\bW_{12}\splus2\bA_4$ admit a
symmetry~$c$ such that the image of~$B$ in $\Cp2(O_c)/c$ is the
maximal trigonal curve
with the set of singularities $2\bA_4$.
However, such sextics form codimension~$1$ subsets in
their equisingular strata, see~\cite{degt.Oka3}.
\endRemark

\Remark\label{conjecture.proof}
One can use a similar parameter count to substantiate the `if'
part of Conjecture~\ref{conjecture}. From the description of the
moduli of sextics, see~\cite{JAG}, it follows that, for any
sextic~$B$ with simple singularities,
$\dim\CM_6(\Sigma_B)=19-\mu(B)$. Let~$B$ be the double covering of
a stable maximal trigonal curve $\B\subset\Sigma_2$ ramified
at~$E$ and a section~$\L$. If $\L$ is transversal to~$\B$, then
$\mu(B)=2\mu(\B)=16$, see Corollary~\ref{mu=8}, and
$\dim\CM_6(\Sigma_B)=3$ is the dimension of the space of sections
of~$\Sigma_2$. Each constraint on~$\L$ (passing though a singular
point of~$\B$, tangency to~$\B$, \etc.) increases $\mu(B)$ by one
while decreasing by one the dimension of the space of sections
(assuming that the constraints are independent). Thus, in all
cases, the dimensions of the moduli spaces coincide, \ie, the deck
translation of the covering is a symmetry of each generic curve in
$\CM_6(\Sigma_B)$.
\endRemark

\widestnumber\key{EO1}
\refstyle{C}

\enddocument